\documentstyle{article}
\newtheorem{defn}{Definition}
\newtheorem{thm}{Theorem}
\newtheorem{bis}{Theorem}
\newtheorem{prop}[thm]{Proposition}
\newtheorem{lem}[thm]{Lemma}
\newtheorem{cor}[thm]{Corollary}

\newenvironment{proof}{\medskip {\bf Proof.}}{\hfill \rule{.5em}{1em}
\\}
\def\Bbb{\bf}

\def\eea{\end{eqnarray*}}
\def\bea{\begin{eqnarray*}}

\begin{document}

\title{On the Scalar Curvature of Einstein Manifolds}

\author{Fabrizio Catanese\\
 Georg-August-Universit\"at G\"ottingen \\
and\\
Claude LeBrun\thanks{Supported 
in part by  NSF grant DMS-9505744.} \\  
SUNY at Stony Brook}

\date{}
\maketitle

\begin{abstract} 
 We show that 
there are high-dimensional smooth
compact manifolds which   admit 
  pairs of Einstein metrics for which the
scalar curvatures have opposite signs. 
These are counter-examples to a conjecture
considered by  Besse \cite[p. 19]{bes}. 
The proof hinges on showing 
that the Barlow surface   has
small deformations with ample  canonical 
line bundle. 
  \end{abstract}

\vfill
\pagebreak

\section{Introduction}

Let $M$ be a smooth compact $n$-manifold. An {\em Einstein
metric} on $M$ is a smooth Riemannian metric $g$ on
$M$ such that 
$$r =\lambda g,$$
where $r $ is the Ricci curvature of $g$ and 
$\lambda$ is any real  constant. Such a
metric  has scalar curvature $s =n\lambda$, so
the   so-called Einstein constant $\lambda$
 and the scalar curvature 
of $g$ have the same sign.  

When $n < 4$, any Einstein metric has constant 
sectional curvature $\lambda /(n-1)$. In dimensions
$2$ and $3$, the sign of $\lambda$ is therefore
determined by the size of the fundamental group, 
and so is really a topological invariant of $M$. 
This motivated Besse \cite[p.19]{bes} to consider  
the conjecture that  no smooth compact $n$-manifold can ever admit
Einstein metrics with different signs of $\lambda$.

The K\"ahler-Einstein case superficially 
seems to support such a conjecture. 
 Indeed, if $J_1$ and $J_2$ are
deformation-equivalent complex structures
on an even dimension smooth manifold $M$, 
then 
there cannot be a K\"ahler-Einstein metric with $\lambda > 0$
compatible with $J_1$ and a K\"ahler-Einstein metric with
$\lambda < 0$ compatible with $J_2$.
When $n=4$, Seiberg-Witten theory \cite{witten,leb}
and the Hitchin-Thorpe inequality \cite{hit} 
allow one to make an even more encouraging 
statement: if a smooth compact 4-manifold   admits a 
 K\"ahler-Einstein metric
with $\lambda \leq 0$, then it cannot admit  {\em any} 
Einstein metric 
for which  $\lambda$ has a different sign. 

 Unfortunately, these pieces of 
 evidence are just red herrings.  
Not only is the conjecture false, but
counter-examples can be constructed as K\"ahler-Einstein
metrics on products of complex surfaces!

\begin{thm} For
any $k \geq 2$, there is a  smooth compact $4k$-manifolds $M$
which  admits both an Einstein metric $g_1$ with $\lambda > 0$
and an Einstein metric $g_2$ with $\lambda < 0$.
Moreover, $M$ may be taken to be the $k$-fold product
of ${\Bbb CP}_2\#8\overline{\Bbb CP}_2$ with itself,
and $g_1$ and $g_2$ may both be taken to be K\"ahler-Einstein
metrics.    
\end{thm}

It follows, of course, that the relevant complex structures
    $J_1$ and $J_2$ on $M$ are   
deformation-inequivalent.

 The    proof begins with the observation that
${\Bbb CP}_2\# 8 \overline{\Bbb CP}_2$
is $h$-cobordant to the   complex surface  $S$ of general type 
discovered by Barlow \cite{bar}.
While 
Barlow's explicit complex structures on $S$ do not have 
ample canonical bundle,  we show in \S\S \ref{first}--\ref{sec}
that suitable small deformations of these  complex structures
do indeed have 
$c_1 < 0$. The theory the complex Monge-Amp\`ere equation 
\cite{aubin,tian,ty,yau}
and standard facts about $h$-cobordisms \cite{smale,wall}
  then imply that  Einstein's equation can be solved  on
$S\times S$   both when  $\lambda >0$ 
and when $\lambda < 0$.

Incidentally, the low-dimensional failure of the $h$-cobordism
theorem is illustrated  by the 
fact that  $S$ and ${\Bbb CP}_2\# 8 \overline{\Bbb CP}_2$ are themselves
non-diffeomorphic \cite{kot,ov}; cf. \cite{leb2}. 
On the other hand, 
Freedman's
work \cite{f} tells us that 
just enough of the $h$-cobordism theorem survives
to conclude that the Einstein manifolds
 $S$ and ${\Bbb CP}_2\# 8 \overline{\Bbb CP}_2$ are {\em homeomorphic}. 
Thus, while the Besse conjecture still stands in dimension 4,
it survives only by virtue of the existence of 
exotic differentiable structures on 4-manifolds.

 So far as we know, these are also the first examples
of smooth compact manifolds which admit 
both K\"ahler metrics of 
positive Ricci curvature and 
K\"ahler metrics of  negative Ricci curvature. 
However,   Lohkamp \cite{loh}, generalizing an earlier result of 
Gao and Yau \cite{gao}, 
 has shown that absolutely every 
smooth  manifold of dimension $n > 2$ admits 
Riemannian metrics of negative
Ricci curvature, so it is only the 
K\"ahler condition that gives this last  observation
 any new interest.  Rethinking the 3-dimensional 
case in light of these last results should  
serve to remind the reader  that the 
 Einstein condition  is immeasurably
 stronger than the mere requirement 
  that the 
Ricci curvature have a fixed sign.

\section{Deformations of the Barlow Surface}
\label{first}

Let $F$ be a quintic surface in ${\Bbb CP}_3$ with 
exactly 20 nodes and no other singularities. 
Each of the  nodes is locally
modeled by ${\Bbb C}^2/{\Bbb Z}_2$, 
since the quadratic 
map ${\Bbb C}^2\to \odot^2{\Bbb C}^2: v\mapsto v\otimes v$
identifies ${\Bbb C}^2/{\Bbb Z}_2$ with the 
nodal surface $\det =0$ in $\odot^2{\Bbb C}^2\cong {\Bbb C}^3$. 
Thus 
$F$  naturally carries the structure of a complex orbifold. 
We will be interested in the case in which
$F$ is a global orbifold:
$$F= Y/{\Bbb Z}_2$$
for some compact complex 2-manifold $Y$
with  a holomorphic involution. Using the 
theory of Hilbert modular surfaces, one can \cite{cat}, for example,
 show that this is true of 
the specific 20-nodal quintic
surface  
\begin{eqnarray} \sum_{j=1}^5 z_j^5&=&\frac{5}{4}
 (\sum_{j=1}^5 z_j^2)( \sum_{j=1}^5 z_j^3) \label{symm}\\
\sum_{j=1}^5 z_j &=& 0 \nonumber ,\end{eqnarray}
 in a hyperplane of ${\Bbb CP}_4$. 
For more details, see \S \ref{sec}.

We will need two basic facts about the surface $Y$.
The first of these is rather obvious:

\begin{lem}
In the above situation, $Y$ has ample canonical line
bundle.  
\end{lem}
\begin{proof}
The canonical line bundle 
bundle of $Y$ is the pull-back of ${\cal O}(1)$
from ${\Bbb CP}_3$ via the tautological map 
$Y\to F$. Since this map is finite-to-one,
it does not collapse any 
curves. The canonical line bundle is therefore 
positive by Nakai's criterion.
\end{proof}
 
The second basic  fact is, by contrast, quite delicate:
\setcounter{bis}{\thethm}

\begin{thm}\label{terza}
In the above  situation, $H^2(Y,TY)=0$. 
\end{thm}
\noindent 
In the interest of clarity, we defer the rather long proof
to \S \ref{sec}.

\medskip

 Now consider the particular quintic $F$ given by 
equation \ref{symm}. The symmetric group ${\cal S}_5$ acts
on $F$, and in particular we have an action of
the dihedral group $D_{10}\subset {\cal S}_5$ 
 of pentagonal isometries, generated by 
$(25)(34)$ and $(12345)$. This action lifts
to $Y$ in such a way that \cite{cat} the cyclic subgroup
${\Bbb Z}_5\subset D_{10}$ generated by 
$(12345)$ acts freely on $Y$ and such that \cite{bar}
the involution
$(25)(34)$ acts with exactly 20 fixed   points.
The so-called Catanese surface $X=Y/{\Bbb Z}_5$ 
is therefore
non-singular, and comes equipped with an
involution $\alpha : X\to X$ with exactly 
4 fixed points. From the above basic facts
we immediately read off

\begin{lem}
The Catanese surface $X$ has ample canonical 
bundle and satisfies $H^2(X,TX)=0$.
\label{cs}
\end{lem}
\begin{proof} These follow from general facts about
unramified covers. For example, 
since $-c_1(Y)$ is represented by a positive 
2-form, it is also, by averaging, represented by 
a positive ${\Bbb Z}_5$-invariant 2-form;
such a form descends to $X$, and represents $-c_1(X)$.
In the same vein, the pull-back maps such as 
$H^2(X,TX)\to H^2(Y,TY)$ are necessarily injective,
since the pull-back of a harmonic representative
is harmonic with respect to an invariant metric. 
\end{proof}

Dividing $X$ by the action of the 
the involution $\alpha$, we obtain 
a surface $X/{\Bbb Z}_2$ whose
only singularities are four nodes. 
The {\em Barlow surface} is by definition
   the minimal resolution  $S$  
of  $X/{\Bbb Z}_2$. One can show \cite{bar}
that  $S$ is a minimal,
simply connected  complex 
surface of general type,
with 
$c_1^2=1$ and $p_g=0$. 

By construction, $S$ contains four 
$(-2)$-curves $B_1, \ldots , B_4$,
arising from the nodes of $X/{\Bbb Z}_2$;
in particular, $S$ does {\em not} have 
ample canonical bundle. On the other hand, 

\begin{lem} The only  $(-2)$-curves in $S$ are 
  $B_1, \ldots , B_4$.
\end{lem}
\begin{proof}
Because  the orbifold $X/{\Bbb Z}_2$ has $K$ ample,
 the underlying variety has  a projective embedding 
in projective space \cite{bail} which lifts to a
pluri-canonical map of $S$ collapsing only the 
$B_j$. 
\end{proof}

Let us  use $B$ to  denote the union of 
the four $(-2)$-curves $B_1, \ldots , B_4$
in $S$, and let $TS (- \log B)$ denote the
sheaf of vector fields on $S$ with trivial
normal component along the curve $B$. 

\begin{lem} $H^2(S, TS (- \log B))=0$.
\label{bs}
\end{lem}
\begin{proof}
 The Serre dual of $H^2(S, TS )(- \log B)$
is $H^0(S, \Omega^1( \log B) \otimes \Omega^2)$.
 The latter  is  \cite[Prop. 3.1]{cat2} exactly the 
${\Bbb Z}_2$-invariant subspace of
$H^0(\tilde{X}, \Omega^1 \otimes \Omega^2)=
H^0({X}, \Omega^1 \otimes \Omega^2)$. 
However, the latter is Serre dual to 
$H^2(X, TX)$, which vanishes by Lemma \ref{cs}. 
\end{proof}

\begin{thm} The Barlow surface $S$ has a
  smooth versal deformation space
of complex dimension 8. A general
point in this space corresponds to a surface 
with ample canonical bundle.
\label{key}
\end{thm}
\begin{proof}
Let $\nu$ denote the normal bundle of 
$B$. We then have an exact sequence
$$0\to  TS(- \log B)\to  TS \to {\cal O}_B(\nu)\to 0$$
of sheaves on $S$. Thus Lemma \ref{bs}
tells us, in particular, that 
$H^2(S,TS)=0$. 
 Kodaira-Spencer theory thus
gives us a smooth versal deformation of 
$S$ with tangent space $H^1(S,TS)$, the dimension of which is
$h^1(S,TS)=-\chi (S, TS) =  (5c_2-7c_1^2)/6= (5\cdot 11 - 7\cdot 1)/6=
8.$ (Note that $H^0(S,TS)=0$ because 
$S$ is of general type.)

Now for each $j=1, \ldots , 4$, define $TS (-\log B_j)$ to be the 
sheaf of holomorphic vector fields with trivial
normal component along the curve $B_j$. We then
have   exact sequences 
$$0\to TS(- \log B)\to TS(-\log B_j)
\to {\cal O}_{
\bar{B}_j}(\nu)\to 0,
$$
where $\bar{B}_j=\cup_{j\neq k}B_k$ 
is the complement of $B_j\subset B$.
Thus 
 the natural maps $H^1(S, TS (-\log B_j))\to
H^1(S,TS)$ have images of codimension 1.
Choose a curve through the base-point of  the Kodaira-Spencer
family, parameterized by an embedding of the unit disk
$\Delta\subset {\Bbb C}$.  
Let ${\cal S}\to \Delta$ be the induced family 
with central fiber $S$. By construction, the Kodaira-Spencer
map of this family at zero does not take values in
$H^1(S,TS (-\log B_j))$ for any $j\in \{1, 2,3 , 4\}$. 
Now   there are only finitely
many classes   $a\in H^2(S, {\Bbb Z})$ such that  $c_1\cdot a=0$
and  $a\cdot a=-2$, since $c_1^2> 0$ and
the intersection form is of Lorentz type. Since
$q(S)=0$, each such $a$ corresponds to a unique 
holomorphic line bundle $L_{a,t}\to S_t$
for each $t\in \Delta$, where $S_t$ is the fiber over $t$.  
We have $h^0(S_t, L_{a,t})\leq 1$ for all $(a,t)$, since
$(-2)$ curves are necessarily isolated, so 
 semi-continuity \cite[\S 10.5]{gr} implies that, for each $a$,
 the set of $t\in \Delta$
where $h^0(S_t,L_{a,t}) = 1$ is closed in the 
analytic Zariski topology. Now suppose there were such 
an $a$
such that $0\in \Delta$ wasn't isolated in
$\{ t\in \Delta ~|~ h^0(S_t,L_{a,t}) = 1\}$.
Then we would have $h^0(S_t,L_{a,t}) = 1$ 
$\forall t\in \Delta$, and the  direct
image theorem \cite[\S 10.5]{gr} 
would assert the existence
of a   section of $L_{t,a}$ which is non-trivial for each $t$ and 
depends holomorphically on $t$. Thus there
would be a family of $(-2)$-curves in
the fibers of ${\cal S}\to \Delta$ depending
smoothly on the base, and in particular
 it would follow that $a=[B_j]$ for 
some $j$. Moreover, we would have a 
family of curves $B_{j,t}\subset S_t$ depending
holomorphically on $t$, so that ${\cal S}$ actually
comes from a deformation of the pair $(S,B_j)$.
But  its Kodaira-Spencer map 
at the base-point would therefore take values in 
$H^1(S, TS (-\log B_j))$, in contradiction to our
assumptions. This shows that there is a non-empty
subset of points in the Kodaira-Spencer space
corresponding to surfaces without $(-2)$-curves;
and by semi-continuity, this non-empty subset 
is the complement of an analytic subset. It follows
\cite{bpv} 
that the general point of the Kodaira-Spencer space 
corresponds to a surface with ample canonical bundle.  
\end{proof}

For the experts, the  above elementary argument may be considerably 
shortened  by
invoking the results of  \cite{bw,cat3}. Indeed, 
from this point of view, Lemma \ref{bs}
asserts  
that there is a smooth submersion
$\mbox{Def }({X/{\Bbb Z}_2})\stackrel{\psi}{\to} {\cal L}_{X/{\Bbb Z}_2}$, 
where ${\cal L}_{X/{\Bbb Z}_2}$ is the local 
deformation space of the four nodes
of ${X/{\Bbb Z}_2}$. Thus by \cite[(1.8)]{cat3}, 
we have a fiber product of smooth manifolds
$$
\begin{array}
{ccc}
\mbox{Def }(S) &\to &{\cal L}_S\\
\downarrow && 
\downarrow_{\beta}\\
\mbox{Def }({X/{\Bbb Z}_2})&\stackrel{\psi}{\to} &{\cal L}_{X/{\Bbb Z}_2}
\end{array}
$$
where $\beta : {\Bbb C}^4\to {\Bbb C}^4$
is given by $(z_j)\mapsto (z_j^2)$. Thus
a general deformation of $S$ has a smooth canonical 
model --- i.e. has ample canonical bundle.

\section{Vanishing Theorems for Nodal Surfaces}
\label{sec}

In this section, we will prove  the key vanishing result  
Theorem \ref{terza}. We will do this 
by  first proving some
general results concerning the following   class of 
nodal surfaces:

\begin{defn}
Let $x\mapsto [a_{ij}(x)]$ be a linear map from 
${\Bbb C}^4$ to the symmetric $n\times n$ matrices,
and let $F\subset {\Bbb CP}_3$ denote the 
 surface of degree $n$ defined by $\det ([a_{ij}(x)])=0$.
Suppose that the only singular points at which
$\mbox{\rm rank} [a_{ij}(x)] \leq n-2$
 are nodes.
Then we shall say that $F$ is a {\em nodal linearly 
symmetric surface} of degree $n$.
\end{defn}

It turns out that there are precisely
 ${n+1\choose 3}$ nodes \cite{cat}, which we shall
denote by  $p_j\in F$, 
$j=1,\ldots,  {n+1\choose 3}$. Let $\tilde{F}$
be the smooth surface obtained by blowing up these nodes,
let $\pi : \tilde{F}\to F$ be the blowing-down map, and set
$A_j=\pi^{-1}(p_j)$. The $A_j$ are thus disjoint $(-2)$-curves,
and we denote their union by $A=\cup_j A_j$. 
Setting ${\cal O}_{\tilde{F}}(H)= \pi^*{\cal O}_F(1)$ and 
  $r=[\frac{n}{2}]$, there is then a divisor $L$ on 
$\tilde{F}$ such that 
$2L=A+ (2r+1-n)H$. It was  shown in 
\cite[Theorem 2.19]{cat} 
that $h^0(\tilde{F}, {\cal O}(rH-L))=n$, and 
the base locus of the system $|rH-L|$ is empty. 
We thus have a holomorphic map
$\varepsilon : \tilde{F}\to {\Bbb CP}_3\times {\Bbb CP}_{n-1}$
given by $|H|\times |rH-L|$, and $\varepsilon$ is an embedding,
since $|H|$ alone suffices away from the nodes and 
$|rH-L|$ is a degree-$1$ embedding on each $A_j$. 
 It was also  shown that, relative to a suitable basis 
$v_j$, $j=1,\ldots , n$ for  $H^0(\tilde{F}, {\cal O}(rH-L))$
that the image of $\varepsilon$ is exactly given by
the equations
$$\sum_{j=1}^n a_{ij}(x)v_j =0, ~i=1,\ldots , n.$$
As was pointed out by L. Ein, this proves the following
useful fact:

\begin{lem} \label{cap}
The blow-up $\tilde{F}$ of such a surface $F$ of degree $n$ 
at its nodes 
is a 
complete intersection of $n$ hypersurfaces in 
${\Bbb CP}_3\times {\Bbb CP}_{n-1}$ of bidegree
$(1,1)$. 
\end{lem}
In particular, if $\tilde{F}$ is smooth,
its   normal bundle in ${\Bbb CP}_3\times {\Bbb CP}_{n-1}$ 
is just the restriction of
$[{\cal O}(1,1)]^{\oplus n}$ to $\tilde{F}$. 
We also get  
the following useful vanishing result:

\begin{lem}
\label{zip} 
If $n > 4$, then
$$ H^1 (\tilde{F}, {\cal O}(a,b)) =0$$ 
provided
that 
\begin{itemize}
\item  $a,b\geq -1$; or
\item $b < 2$, $a < n-2$; or 
\item $-1\leq a \leq n-3$; or
\item $-1\leq b\leq 1$.
\end{itemize}
Moreover, 
$$H^0({\Bbb CP}_3\times {\Bbb CP}_{n-1},
{\cal O}(a,b))\longrightarrow H^0(\tilde{F},
{\cal O}(a,b))$$
is surjective provided
that 
\begin{itemize}
\item $a,b\geq 0$; or
\item $b < 3$, $a < n-1$; or 
\item $0\leq a \leq n-2$; or
\item $0\leq b\leq 2$.
\end{itemize}
\end{lem}
\begin{proof}
Let $E=[{\cal O}_{{\Bbb CP}_3\times 
{\Bbb CP}_{n-1}}(-1,-1)]^{\oplus n}$.
Because $\tilde{F}$ is a complete
intersection, we have the 
Koszul complex
$$
0\to \Lambda^nE\to \cdots \to \Lambda^2E\to E\to 
{\cal O}_{{\Bbb CP}_3\times 
{\Bbb CP}_{n-1}}\to 
{\cal O}_{\tilde{F}}\to 0,$$
and we may tensor this with ${\cal O}(a,b)$ to obtain an
exact complex 
$$0\to {\cal E}^{-n}\to \cdots \to {\cal E}^{-2}\to {\cal E}^{-1} 
\to {\cal O}(a,b)\to {\cal O}_{\tilde{F}}(a,b)\to 0.$$
Since a sum of line bundles on a projective space can only have 
non-trivial cohomology in the top or bottom dimension, 
 the K\"unneth formula tells us that this complex always satisfies
 $H^i({\cal E}^{-j})=0$ when $i\not\in\{ 0,3,n-1 , n+2\}$.
By breaking the complex into  short exact sequences, it thus follows 
that $H^1(\tilde{F}, {\cal O}(a,b))=0$
whenever  
 $H^3({\cal E}^{-2})=H^{n-1}({\cal E}^{-(n-2)})=0$.
 But ${\cal E}^{-2}$ is a direct sum of copies of 
${\cal O}(a-2,b-2)$, so that $H^3({\cal E}^{-2}) =0$
provided that $a > -2$ or $b < 2$. Similarly, 
since ${\cal E}^{-(n-2)}$ is a direct sum of copies of 
${\cal O}(a-(n-2),b-(n-2))$, we have 
$H^{n-1}({\cal E}^{-(n-2)})=0$ provided 
that $a < n-2$ or $b > -2$. 

The same splicing argument
shows that the surjectivity of  
$$H^0({\Bbb CP}_3\times {\Bbb CP}_{n-1},
{\cal O}(a,b))\longrightarrow H^0(\tilde{F},
{\cal O}(a,b))$$
 follows whenever  
$H^3({\cal E}^{-3})=H^{n-1}({\cal E}^{-(n-1)})=0$.
Since ${\cal E}^{-3}$ and ${\cal E}^{-(n-1)}$
are respectively   sums of copies of 
${\cal O}(a-3,b-3)$
and ${\cal O}(a-(n-1),b-(n-1))$,
one has $a > -1$ or $b < 3 ~\Rightarrow H^3({\cal E}^{-3})=0$,
and $a < n-2$ or $b > -1  ~\Rightarrow H^{n-1}({\cal E}^{-(n-1)})=0$.
\end{proof}

This, incidentally, yields a sharpening 
 of  \cite[Proposition 2.25]{cat}:

\begin{cor}
Let $p_j\in {\Bbb CP}_3$ be the node points of $F$, as before,
and for each $j=1, \ldots , {n+1\choose 3}$, consider the 
condition of passing through $p_j$
as a linear constraint on the linear system 
of surfaces of degree $(n-1)$ in ${\Bbb CP}_3$.
Then these ${n+1\choose 3}$ conditions are linearly independent. 
\end{cor}
\begin{proof}
The restriction map 
$$H^0({\Bbb CP}_3\times {\Bbb CP}_{n-1}, {\cal O}(0,2))\to
H^0(\tilde{F}, {\cal O}(0,2))$$
is surjective by Lemma \ref{zip},
and it is injective by the same reasoning. However, 
${\cal O}_{\tilde{F}}((n-1)H-A)={\cal O}_{\tilde{F}}(0,2)$
by construction.
Thus
$$h^0(\tilde{F}, {\cal O}((n-1) H-A))=
h^0({\Bbb CP}_3\times {\Bbb CP}_{n-1}, {\cal O}(0,2))=
{n+1\choose 2}.$$
But the restrict-and-pull-back 
 map $H^0({\Bbb CP}_3, {\cal O}(n-1))\to 
H^0(\tilde{F}, {\cal O}((n-1) H))$ is also an isomorphism,
so $H^0(\tilde{F}, {\cal O}((n-1) H-A))$ can be
identified with the subspace of 
$H^0({\Bbb CP}_3, {\cal O}(n-1))$ consisting of 
sections which vanish at the $p_j$. Since
$${
n+1 \choose 2}= 
{
n+2\choose 3}-
{n+1\choose 3}=
h^0({\Bbb CP}_3, {\cal O}(n-1)) - \# \{ p_j\}, $$
the claim therefore follows. 
\end{proof}

We now specialize to the $n=5$ case. Let 
$F$ be a linearly symmetric quintic whose only 
singularities are its 20 nodes;  an example
is \cite{cat} given by \ref{symm}.  Thus $\tilde{F}$ 
is smooth. Moreover, $A=2L$, so there is a double 
branched cover $\tilde{Y}$ ramified exactly over the 
exceptional curves $A_j$. The inverse image of
each $A_j$ is an $(-1)$-curve $E_j$, and we
may blow these 20 disjoint $(-1)$-curves down to 
obtain a smooth surface $Y$, which comes equipped with
a natural involution $\Phi: Y\to Y$
corresponding to the
sheet-interchanging map $\tilde{\Phi}: \tilde{Y}\to 
\tilde{Y}$  of $\tilde{Y}\to \tilde{F}$. 
We thus have $F=Y/{\Bbb Z}_2$, realizing $F$ as a 
global orbifold, exactly as claimed in \S \ref{first}. 

We are now prepared to prove the key vanishing result
used in the previous 
section:

\begin{bis}
Let $Y$ be as above. Then $H^2(Y,TY)=0$.
\end{bis}
\begin{proof}
The Serre dual of $H^2(Y,TY)$ is $H^0(Y, \Omega^1\otimes \Omega^2)$,
and the latter can be identified with 
$H^0(\tilde{Y}, \Omega^1\otimes \Omega^2)$
as a consequence of  Hartogs' theorem. Let us express 
the latter in terms of the $(\pm 1)$-eigenspaces
of the action of $\tilde{\Phi}$:
$$H^0(\tilde{Y}, \Omega^1\otimes \Omega^2)
=H^0(\tilde{Y}, \Omega^1\otimes \Omega^2)^+\oplus
H^0(\tilde{Y}, \Omega^1\otimes \Omega^2)^-.$$ 
By standard formulae \cite[Prop. 3.1]{cat2},
these eigenspaces have direct interpretations on $\tilde{F}$: 
\bea
 H^0(\tilde{Y}, \Omega^1\otimes \Omega^2)^+&=&
 H^0(\tilde{F}, \Omega^1(\log A)\otimes \Omega^2) 
\\
H^0(\tilde{Y}, \Omega^1\otimes \Omega^2)^-&=&
H^0(\tilde{F}, \Omega^1 \otimes \Omega^2\otimes L). 
\eea
In fact, however, by \cite[Proposition 1.6]{cat3}, 
the first statement has an ostensible 
  strengthening: 
$$H^0(\tilde{Y}, \Omega^1\otimes \Omega^2)^+=
H^0(\tilde{F}, \Omega^1 \otimes \pi^*\omega_F) =
H^0(\tilde{F}, \Omega^1 \otimes \Omega^2) .$$
This strengthening is Serre dual to 
the equality 
$$H^2(\tilde{F}, T\tilde{F}(-\log A))=
H^2(\tilde{F}, T\tilde{F}),$$
which follows directly \cite[Corollary 1.3]{bw}
from the observation that ${F}$ has deformations which
smooth all its nodes.

Let us now compute $H^0(\tilde{Y}, \Omega^1\otimes \Omega^2)^+=
H^0(\tilde{F}, \Omega^1 \otimes \Omega^2\otimes L)$.
By Lemma \ref{cap}, $\tilde{F}$ is the complete
intersection of five hypersurfaces in ${\Bbb CP}_3\times 
{\Bbb CP}_4$ of bidegree $(1,1)$. The canonical line
bundle of $\tilde{F}$ is therefore
$$\Omega^2_{\tilde{F}}= {\cal O}_{\tilde{F}}(-4,-5)\otimes 
{\cal O}_{\tilde{F}}(5,5)= {\cal O}_{\tilde{F}}(1,0),$$
whereas 
$${\cal O}(L)= {\cal O}(2H) \otimes
[{\cal O}(2H-L)]^*  = {\cal O}_{\tilde{F}}(2,-1).$$
Thus our objective is to show that 
$H^0( \Omega^1_{\tilde{F}} (3,-1))$ vanishes. 
Because $H^1(\tilde{F}, {\cal O}(2,-2))=0$
by  Lemma \ref{zip}, tensoring the conormal 
bundle sequence 
$$0\to [{\cal O}_{\tilde{F}} (-1,-1)]^{\oplus 5} \to 
 \hat{\Omega}^1
\to \Omega^1_{\tilde{F}}\to 0$$
with ${\cal O}(3,-1)$
shows that it is enough to ascertain the
vanishing of 
$H^0(\tilde{F},  \hat{\Omega}^1 (3,-1) )$,
where   $\hat{\Omega}^1=\Omega^1_{{\Bbb P} \times{\Bbb P}}
\otimes {\cal O}_{\tilde{F}}$ 
is the restriction of the cotangent bundle of 
${\Bbb CP}_3\times{\Bbb CP}_4$ to $\tilde{F}$. 
However,   tensoring the
Euler exact sequence
$$0\to \hat{\Omega}^1 \to 
[{\cal O}(-1,0)]^{\oplus 4}\oplus [{\cal O}(0,-1)]^{\oplus 5}\to
{\cal O}\oplus {\cal O}\to 0$$
with ${\cal O}_{\tilde{F}}(3,-1)$ tells us that 
$$  H^0(\tilde{F}, \hat{\Omega}^1 (3,-1) 
 )\subset 
[H^0(\tilde{F}, {\cal O}(3,-2))]^{\oplus 4}\oplus 
[H^0(\tilde{F}, {\cal O}(2,-1))]^{\oplus 5} .$$
On the other hand,  $H^0({\Bbb CP}_3\times {\Bbb CP}_4, 
{\cal O}(3,-2))\to H^0(\tilde{F}, 
{\cal O}(3,-2))$ and $H^0({\Bbb CP}_3\times {\Bbb CP}_4, 
{\cal O}(2,-1))\to H^0(\tilde{F}, 
{\cal O}(2,-1))$ are both surjective by Lemma
\ref{zip}; and since the relevant line bundles on 
${\Bbb CP}_3\times {\Bbb CP}_4$ obviously do not
admit non-trivial sections, we conclude that 
$H^0(\tilde{Y}, \Omega^1\otimes \Omega^2)^-=0$.

To finish the proof, we now compute 
$H^0(\tilde{Y}, \Omega^1\otimes \Omega^2)^+=
H^0(\tilde{F}, \Omega^1 \otimes \Omega^2)$
by the same method. Indeed, 
$\Omega^1_{\tilde{F}}\otimes \Omega^2_{\tilde{F}}=
 \Omega^1_{\tilde{F}}(1,0)$, and the conormal sequence 
tells us that 
$$H^0(\tilde{F},  \hat{\Omega}^1(1,0)) \to 
H^0(\tilde{F},  \Omega^1(1,0))\to
[H^1 (\tilde{F}, {\cal O}(0,-1))]^{\oplus 5}  $$
is exact. Since  $H^1 (\tilde{F}, {\cal O}(0,-1))=0$ by
Lemma \ref{zip}, we need only verify that 
$H^0(\tilde{F},  \hat{\Omega}^1(1,0))$ vanishes, too.
Now the Euler sequence tells us that 
$$  
H^0(\tilde{F},  \hat{\Omega}^1(1,0))= \ker  
\left( 
\begin{array}{c}
 [ H^0(\tilde{F},  {\cal O})]^{\oplus 4} \cr 
 \oplus  \cr 
 [ H^0 
( \tilde{F},  {\cal O}(1,-1)) ]^{\oplus 5}
\end{array} 
\right) \to 
\left( \begin{array}{c}
H^0 (\tilde{F},  
{\cal O}(1,0)) \\ 
\oplus \\ 
 H^0 (\tilde{F},  {\cal O}(1,0))
\end{array}
\right)   .
$$
However, $H^0(\tilde{F},  {\cal O}(1,-1))=0$
by Lemma \ref{zip}, and $[H^0(\tilde{F},  {\cal O})]^{\oplus 4}\to
H^0(\tilde{F},  {\cal O}(1,0))$ is manifestly injective. 
Thus 
$H^0(\tilde{Y}, \Omega^1\otimes \Omega^2)^+=0$,
and the result follows.
\end{proof}

It is worth noting that
the above result depends heavily on the assumption that 
$F$ is quintic. Indeed, one can use  much the same argument
to show that the 
the  vanishing assertion definitely  fails iff $n \gg 5$.

\section{Einstein Manifolds and Smooth Topology}

We begin by recalling that two smooth compact
oriented $n$-manifolds are said to be 
$h$-cobordant \cite{smale} if there is a compact oriented
$(n+1)$-manifold $W$ whose boundary is $M\cup \bar{N}$,
and such that the inclusion of either boundary component
is a homotopy equivalence; here $\bar{N}$ denotes the 
orientation-reversed version of $N$. One then says
that $W$ is an $h$-cobordism from $M$ to $N$. 
This defines an equivalence relation on the 
set of all $n$-manifolds, since gluing two $h$-cobordisms 
end-to-end yields a new $h$-cobordism.  
This equivalence relation is also compatible with 
Cartesian products:

\begin{lem}
Let $M$ and $N$ be $h$-cobordant   manifolds.
 Then their Cartesian self-products 
$M\times M$ and $N\times N$ are also 
$h$-cobordant. Moreover, the 
 iterated self-products
$M^{(k)}=M\times \cdots \times M$ and 
$N^{(k)}=N\times \cdots \times N$ are
$h$-cobordant for any $k\geq 1$. 
\end{lem}
\begin{proof}
If $W$ is a an $h$-cobordism between $M$ and $N$, 
then $M\times W$ is an $h$-cobordism between $M\times M$ and 
$M\times N$, whereas $W\times N$  
is an $h$-cobordism between $M\times N$ and 
$N\times N$; gluing these end-to-end then
gives the desired $h$-cobordism.

The iterated case follows   similarly, using
induction. Namely, if $W_{k-1}$
is an $h$-cobordism from $M^{(k-1)}$ to
$N^{(k-1)}$, then $M\times W_{k-1}$
is an $h$-cobordism between $M^{(k-1)}$  and 
 $M\times N^{(k-1)}$,
whereas $W\times N^{(k-1)}$ is an 
$h$-cobordism between $M\times N^{(k-1)}$ 
and $N^{(k)}$. 
\end{proof}

This observation allows one to prove the following useful fact: 

\begin{prop}
Let $M$ and $N$ be two  smooth simply connected 
compact 4-manifolds
with equal Euler characteristics
$\chi (M)=\chi (N)$ and   signatures
$\tau (M)=\tau (N)\not\equiv 0 \bmod 16$.
Then $M\times M$ is diffeomorphic to $N\times N$. Moreover,
the   iterated Cartesian self-products $M^{(k)}$
and $N^{(k)}$ are diffeomorphic 
for any $k > 1$. 
\end{prop}
 \begin{proof}
By Rochlin's theorem \cite{ah}, $M$ and $N$ are not spin, and so have
odd intersection form. Since $b_2(M)=b_2(N)$
and $\tau (M)=\tau (N)$, the intersection forms of 
$M$ and $N$ are therefore isomorphic by the 
Minkowksi-Hasse classification \cite{hm} and Donaldson's
thesis \cite{don}. 
A  theorem of Wall \cite{wall}    
therefore tells us that $M$ and $N$ are
$h$-cobordant. By the previous lemma,
it follows that 
of $M^{(k)}$ and $N^{(k)}$ are $h$-cobordant for any 
$k$. But for any $k > 1$, these   
are simply connected compact smooth manifolds
of dimension $> 5$, so Smale's $h$-cobordism theorem
  asserts that they are diffeomorphic, as claimed.
\end{proof}

\begin{cor}
Let $S$ denote the Barlow surface, and let $R$
denote the rational complex surface ${\Bbb CP}_2\#
8 \overline{\Bbb CP}_2$.
Then $S^{(k)}$ is diffeomorphic to $R^{(k)}$
for all $k > 1$. In particular, $S\times S$ is diffeomorphic 
to $R\times R$.
\end{cor}
\begin{proof}
Both $R$ are $S$ are simply connected complex surfaces with 
$c_1^2=1$ and $p_g=0$. Thus both have $\chi =11$ and
$\tau = -7 \not\equiv 0 \bmod 16$. The claim thus follows
from the preceding result.
\end{proof}

\begin{thm}
Let $R={\Bbb CP}_2\#
8 \overline{\Bbb CP}_2$. Then the
8-manifold $R\times R$ admits Einstein
metrics of both positive and negative scalar
curvatures. Similarly, the $4k$-manifold
$R^{(k)}$ admits Einstein
metrics of both positive and negative scalar
curvatures for all $k > 1$. Moreover, the
relevant 
Einstein metrics can be constructed so
as to be K\"ahler with respect to suitable 
complex structures. 
\end{thm}
\begin{proof}
By Theorem \ref{key}, the Barlow surface $(S,J_0)$ has 
deformations with ample canonical line
bundle. 
If $J_t$ is any complex structure on $S$ with this
property,   the work of Aubin and Yau \cite{aubin,yau}
tells us that there is a unique Einstein metric $g_S$
on $S$ which is K\"ahler with respect to $J_t$ and
has scalar curvature $-1$. On the other hand,
if $\tilde{J}$ is a complex  structure
on $R$ corresponding to a  blow-up of
${\Bbb CP}_2$ at 8 points in general position,
then Tian \cite{tian}, building on his joint
work with Yau \cite{ty}, has shown that 
$R$ admits an Einstein metric $g_R$ of scalar
curvature $+1$ which is K\"ahler with respect
to $\tilde{J}$. 

Now because the Ricci tensor of any Riemannian
product is the direct  sum of the Ricci tensors
of the factors, any Riemannian product of two 
Einstein manifolds with the same value of $\lambda$
is again Einstein. In particular, Cartesian self-products
of Einstein manifolds are always Einstein, and, moreover,  
Cartesian self-products of K\"ahler-Einstein manifolds
are always K\"ahler-Einstein. Thus $R^{(k)}$
admits K\"ahler-Einstein metrics with 
positive scalar curvature, and $S^{(k)}$
admits K\"ahler-Einstein metrics with
negative scalar curvature. But by the above Corollary,
these manifolds are diffeomorphic, and the relevant
Einstein metrics may therefore both be 
considered as living on $R^{(k)}$. 
\end{proof}

\bigskip

\noindent
{\bf Acknowledgements.}
The first author would like to 
thank L. Ein for pointing out that $\tilde{F}$
is a complete intersection. He would also like to 
to thank P. Supino for his queries regarding Theorem  
\ref{terza}.

\vfill

\noindent
{\bf Authors' addresses:}

\bigskip

\noindent 
Prof. Fabrizio Catanese\\
Mathematisches Institut\\
Bunsenstra{\ss}e 3/5\\
 37073 G\"ottingen, Germany

\bigskip

\noindent
Prof. Claude LeBrun 
\\ Department of Mathematics\\
SUNY at Stony
 Brook\\
Stony Brook, NY 11794-3651 USA

\end{document}